\documentclass[]{amsart}
\usepackage[utf8]{inputenc}
\usepackage[mathscr]{eucal} 
\usepackage{stmaryrd} 
\usepackage{amsmath,amsthm,amsfonts,amssymb,amscd} 
\usepackage[dvipsnames]{xcolor} 
\usepackage{xspace} 
\usepackage{fancyhdr} 
\usepackage{graphicx} 
\usepackage{listings} 
\usepackage[]{hyperref} 
\usepackage{enumitem} 
\usepackage{relsize} 
\usepackage{tikz-cd} 
\usepackage{mdwlist} 
\usepackage{multicol} 
\usepackage{float} 
\usepackage{adjustbox} 
\usepackage{tikz} 
\usepackage[noabbrev,nameinlink]{cleveref} 
\Crefformat{section}{#2\S#1#3}
\usepackage{setspace} 
\usepackage{bbm} 
\usetikzlibrary{matrix, calc, arrows}

\hypersetup{
    colorlinks=true, 
    linktoc=all,     
    linkcolor=Brown,citecolor=Brown,urlcolor=MidnightBlue,  
}

\setlist[enumerate]{label=(\alph*)} 
\usetikzlibrary{graphs,decorations.pathmorphing,decorations.markings}
\tikzcdset{scale cd/.style={every label/.append style={scale=#1},
        cells={nodes={scale=#1}}}}

\newcommand{\Hom}{\operatorname{Hom}}

\newcommand{\End}{\operatorname{End}}

\newcommand{\id}{\operatorname{id}}

\newcommand{\stab}{\operatorname{stab}}

\newcommand{\Spc}{\operatorname{Spc}}

\newcommand{\pt}{\operatorname{pt}}

\newcommand{\on}[1]{\operatorname{#1}}

\newcommand{\calH}{\mathcal{H}}

\newcommand{\catC}{\mathscr{C}}

\newcommand{\catF}{\mathscr{F}}
\newcommand{\catG}{\mathscr{G}}

\newcommand{\catI}{\mathscr{I}}
\newcommand{\catJ}{\mathscr{J}}
\newcommand{\catK}{\mathscr{K}}
\newcommand{\catL}{\mathscr{L}}

\newcommand{\catP}{\mathscr{P}}

\newcommand{\bbF}{\mathbb{F}}

\newcommand{\bbN}{\mathbb{N}}

\newcommand{\bbZ}{\mathbb{Z}}
\newcommand{\bbone}{\mathbbm{1}}

\newcommand{\up}{{\color{darkblue}\uparrow}}
\newcommand{\down}{{\color{darkblue}\downarrow}}

\def\dot{\color{darkblue}{\color{white}\bullet}\!\!\!\circ}

\definecolor{darkblue}{HTML}{111199}
\definecolor{darkgreen}{HTML}{336633}
\definecolor{darkred}{HTML}{993333}
\definecolor{darkpurple}{HTML}{995599}

\newcommand{\Sp}{\on{Sp}}

\newcommand{\supp}{\on{supp}}
\newcommand{\sbull}{{\scriptscriptstyle\bullet}}

\newcommand{\Ext}{\on{Ext}}
\newcommand{\Spec}{\on{Spec}}
\newcommand{\cp}{^{cp}}

\newcommand{\Proj}{\on{Proj}}
\newcommand{\Heis}{\calH eis}

\newtheorem{theorem}{Theorem}[section]
\newtheorem{lemma}[theorem]{Lemma}
\newtheorem{proposition}[theorem]{Proposition}

\newtheorem*{theorem*}{Theorem}

\newtheorem*{question}{Question}

\newtheorem{innercustomthm}{Theorem}
\newenvironment{customthm}[1]
  {\renewcommand\theinnercustomthm{#1}\innercustomthm}
  {\endinnercustomthm}

\theoremstyle{remark}
\newtheorem{remark}[theorem]{Remark}

\theoremstyle{definition}
\newtheorem{definition}[theorem]{Definition}

\newtheorem*{definition*}{Definition}

\AddToHook{env/corollary/begin}{\crefalias{theorem}{corollary}}
\AddToHook{env/proposition/begin}{\crefalias{theorem}{proposition}}
\AddToHook{env/lemma/begin}{\crefalias{theorem}{lemma}}
\AddToHook{env/definition/begin}{\crefalias{theorem}{definition}}
\AddToHook{env/example/begin}{\crefalias{theorem}{example}}
\AddToHook{env/remark/begin}{\crefalias{theorem}{remark}}
\AddToHook{env/observation/begin}{\crefalias{theorem}{observation}}

\makeatletter
\newcommand{\xhookdoubleheadrightarrow}[2][]{%
  \lhook\joinrel
  \ext@arrow 0359\rightarrowfill@ {#1}{#2}%
  \mathrel{\mspace{-15mu}}\rightarrow
}
\makeatother

\begin{document}
    \title{A semisimple subcategory of Khovanov's Heisenberg category}

    \author[S.~K.~Miller]{Sam K. Miller}
    \address{S.~K.~Miller,
    Department of Mathematics, University of Georgia, Athens GA 30602, United States of America} 
    \email{sam.miller@uga.edu} 

    \begin{abstract}
        We show the existence of a semisimple replete subcategory of Khovanov's Heisenberg category that retains the isomorphism data of objects for the full category. This leads to a noncommutative tensor-triangular geometric example of a monoidal triangulated category whose Balmer spectrum satisfies the tensor product property but which contains one-sided thick tensor-ideals that are not two-sided, and whose standard support varieties fail to classify one-sided thick $\otimes$-ideals.
    \end{abstract}

    \maketitle
    \section*{Introduction}

    Throughout the paper, $\catK$ denotes an essentially small monoidal-triangular category. That is, $\catK$ is a triangulated category with monoidal structure that is exact in each variable, but we need not assume $\catK$ is symmetric. As an appetizer, we'll begin with a noncommutative tensor-triangular geometric theorem.

    \begin{theorem*}{\cite[Theorem 4.1.1]{NVY22b}}
        Let $\catK$ be an essentially small monoidal triangulated category. Then every prime thick $\otimes$-ideal of $\catK$ is completely prime if and only if the \emph{tensor product property} holds for $\catK$, that is, for all $x, y\in \catK$, \[\supp(x) \cap \supp(y) = \supp(x \otimes y).\] Moreover, these equivalent conditions hold if $\catK$ is right duo, i.e., all right thick $\otimes$-ideals are two-sided. (\cite[Theorem 3.2.1]{NVY22b})
    \end{theorem*}

    This so-called \textit{tensor-product property} in support theory has recently emerged as a question of interest. We previously showed in \cite{Mil25} that \textit{complete prime spectrum} $\Spc\cp(\catK)$ is universal amongst support theories for $\catK$ satisfying the tensor product property, even if $\Spc(\catK) = \emptyset$. On the cohomological end, Bergh--Plavnik--Witherspoon showed in \cite{BPW25} that for any finite tensor category $\catC$, the cohomological support variety $\Proj(\Ext^\sbull_\catC(\bbone,\bbone))$ satisfies the tensor product property if and only if it does for pairs of indecomposable periodic objects. They conjectured the tensor product property holds if $\catC$ is braided and verified this in a number of examples. Their conjecture is implied by a stronger conjecture of Nakano--Vashaw--Yakimov, spelled out in detail in \cite{NVY25}. One prediction of the authors is that for any finite tensor category $\catC$ and its stable module category $\catK = \stab(\catC)$, its cohomological support variety and the noncommutative Balmer spectrum are one and the same.

    The noncommutative Balmer spectrum is in general non-functorial, which can make explicit computation difficult in practice. On the other hand the complete prime spectrum is functorial \cite{Mil25}, thus if one knows that given a monoidal triangulated category $\catK$ that every prime is completely prime, one has more tools at their disposal. For instance, we previously extended a result of Balmer \cite{Bal18}: given a conservative functor of monoidal triangulated categories $F\colon\catK \to \catL$, the induced map on complete prime spectra $F^*\colon \Spc\cp(\catL) \to \Spc\cp(\catK)$ is surjective when $\catL$ is duo and $\catK$ is rigid. Given that duoness appears to be a necessary assumption, it is pertinent to understand this property. One natural question to ask is the following:

    \begin{question}
        Does the converse to \cite[Theorem 4.2.1]{NVY22b} hold? That is, if the tensor product property holds for $\catK$, is $\catK$ left or right duo?
    \end{question}

    In this paper, we give an example demonstrating the answer is \textbf{no}, as expected. Though one should not use ring theory as a motivating example in tensor-triangular geometry, there are plenty of noncommutative rings for which all primes are completely prime, but have one-sided ideals that are not two-sided, for instance, matrix rings. To that point, we remark that Solberg--Vashaw--Witherspoon recently initiated a program to understand the noncommutative tensor-triangular geometry of bimodule categories for selfinjective algebras \cite{SVW25}; as bimodules can be viewed as categorifications of homomorphisms, we expect this setting may produce more examples of this behavior.

    \subsection*{Entering the diagrammatic world} Our example arises from a territory which tensor-triangular geometry seems to have not yet encroached: diagrammatic categories. We examine the Heisenberg category $\Heis_k(\delta)$ defined over a ring $\bbF$, first constructed by Khovanov \cite{Kho14} for $k = - 1$ and $\delta = 0$, then quantized by Licata--Savage \cite{LS13, LS15}, and extended to arbitrary central charge by Brundan and Mackaay--Savage \cite{Bru18, MS18}. Khovanov conjectured that the Karoubi envelope (idempotent completion) of $\Heis := \Heis_{-1}(0)$ categorifies the Heisenberg algebra; this was proven by Brundan--Webster--Savage \cite{BSW23}. The Heisenberg category also is closely related to Kac-Moody categorification \cite{BSW20, BSW25}, and at central charge 0 recovers the affine oriented Brauer category \cite{BCNR17, Bru18}. Moreover, $\Heis := \Heis_{-1}(0)$ categorifies the integral \textit{Weyl algebra} on one variable $A_1$; we note this in \Cref{prop:isoclasses}.

    The tensor-triangular geometry of diagrammatic categories appears at present relatively unexplored. We warn the reader that the term \textit{tensor ideal} in the context of tensor categories has a different meaning than the term \textit{thick tensor ideal} in the context of tensor-triangulated categories. In the context of monoidal categories there is an additional weaker notion of a \textit{thick ideal}; computations of these have been performed for Deligne categories (which are diagrammatic), see for instance \cite{CH17,Cou18}.

    The Heisenberg category appears difficult to do tensor-triangular geometry with naively, as it not naturally triangulated, is highly nonabelian and has large morphism spaces. We remedy this situation by removing almost all of the interesting behavior of the category and instead showing the existence of a replete semisimple subcategory $\Heis_s \subset \Heis$ that encodes the isomorphism data on objects; this is \Cref{thm:semisimplesubcat} and \Cref{prop:isoclasses}.

    \begin{customthm}{A}
        There exists a replete semisimple monoidal subcategory $\Heis_{s} \subseteq \Heis $ for which the Heisenberg isomorphism $\down \otimes \up\cong \up \otimes \down\oplus \bbone$ holds in $\Heis_{s}$. The simple objects are given precisely by $\up^{\otimes i}\otimes \down^{\otimes j}$ for nonnegative integers $i,j \geq 0$.

        We have an isomorphism $K_0(\Heis_s) \cong A_1$, and if $\bbF$ is a field of characteristic 0, $K_0(\Heis_s)\cong K_0(\Heis)$. In particular, given two objects $X,Y \in \Heis$, $X \cong Y$ in $\Heis$ if and only if $X \cong Y$ in $\Heis_s$.
    \end{customthm}

    The proof is fairly computational, we restrict down to morphisms in $\Heis$ satisfying certain properties, and show the resulting subcategory is well-defined. We remark that the problem of determining whether a $\bbZ_+$ ring has a semisimple categorification is rather difficult in general; see \cite[Remark 4.10.8]{EGNO15} for discussion. It is unclear if such a semisimple subcategory exists in other central charges or quantizations of the Heisenberg category, but we expect this to not be the case.

    From here, we deduce that the category $\catK := D^b(\Heis_s)$ satisfies the tensor product property, but is not duo, by explicitly computing all of its left, right, and two-sided thick $\otimes$-ideals; this is \Cref{thm:classificationofideals}.

    \begin{customthm}{B}
        Let $\catK := D^b(\Heis_s)$. There are 2 two-sided thick $\otimes$-ideals of $\catK$, $0$ and $\catK$, but the thick left $\otimes$-ideals are given by \[\catK \supset \langle \down\rangle_l \supset \langle \down\otimes\down\rangle_l \supset \cdots\supset 0,\] and the thick right $\otimes$-ideals are given by \[\catK \supset \langle \up \rangle_r \supset \langle \up\otimes\up \rangle_r \supset \cdots\supset 0.\]

        In particular, $\Spc(\catK) = \{\ast\}$, and for all nonzero $x \in \catK$, $\supp(x) = \Spc(\catK)$. Therefore, $\catK$ satisfies the tensor product property but is not left or right duo.
    \end{customthm}

    The fact that the one-sided thick $\otimes$-ideals of $\catK$ form a totally ordered set provides a rare occurrence: the one-sided ideals are not classified by either the noncommutative Balmer spectrum $\Spc(\catK)$ nor the cohomological spectrum $\Spec^h(\End_{\catK}^\sbull(\bbone))$ (which is to be expected), but are classified by the \textit{universal quasi-support datum} $\Sp(\catK)$, which as a set consists of all proper thick one-sided $\otimes$-ideals. We discuss this in \Cref{rmk:classification}.

    The paper is organized as follows: preliminaries are reviewed in Section 1, and the main results are proved in Section 2.

    \subsection*{Acknowledgments} We thank Jon Brundan for his comments and allowing us to use
    his TikZ code. We also thank Math.SE user Qiaochu Yuan for answering a related question on Math.StackExchange \cite{Yua25} which gave inspiration for this paper. The author is partially supported by an AMS-Simons travel grant.

    \section{Preliminaries}

    \subsection{Noncommutative tensor-triangular geometry}
    We quickly lay out the essentials of \textit{noncommutative tensor-triangular geometry}. Nakano--Vashaw--Yakimov laid out the theory in \cite{NVY22}, extending the ideas of Balmer's seminal paper \cite{Bal05} to the setting of non-symmetric tensor products. For an overview of tensor-triangular geometry, see \cite{Bal10b}.

    Recall $\catK$ denotes an essentially small monoidal triangulated category. A \textit{thick subcategory} of $\catK$ is a full triangulated subcategory of $\catK$ closed under direct summands. A thick subcategory $\catI \subseteq \catK$ is a \textit{thick $\otimes$-ideal} if in addition, for all $x \in \catI$ and $y \in \catK$, $x \otimes y \in \catK$ and $y \otimes x \in \catK$. Given an element $x \in \catK$, we write $\langle x\rangle$ (resp.\ $\langle x\rangle_l$, $\langle x\rangle_r$) to denote the thick $\otimes$-ideal (resp.\ thick left $\otimes$-ideal, thick right $\otimes$-ideal) generated by $x$.

    We say $\catK$ is \textit{left (resp. right) duo} if all thick left $\otimes$-ideals (resp. thick right $\otimes$-ideals) are two-sided. If $\catK$ is both left and right duo, we say it is duo.

    \begin{definition}{\cite{NVY22}}
        A proper thick $\otimes$-ideal $\catP \subset \catK$ is \textit{prime} if for all thick $\otimes$-ideals $\catI, \catJ$ of $\catK$, $\catI \otimes \catJ \subseteq \catP$ implies $\catI \subseteq \catP$ or $\catJ \subseteq \catP$. Equivalently $\catP$ is prime if for all $x, y \in\catK$, $x \otimes \catK \otimes y \subseteq \catP$ implies $x \in \catK$ or $y \in \catK$. Finally, $\catP$ is \textit{completely prime} if for all $x, y \in\catK$, $x \otimes y \in \catP$ implies $x\in\catP$ or $y \in\catP$. If $\catK$ is symmetric (or braided), every prime ideal is completely prime.

        The \textit{noncommutative Balmer spectrum} $\Spc(\catK)$ of an essentially small monoidal-triangulated category $\catK$ is, as a set, the collection of prime thick $\otimes$-ideals of $\catK$ equipped with the topology defined by taking the supports of objects \[\supp(x) := \{\catP \in \Spc(\catK) \mid x \not\in \catP\}\] as a basis of closed subsets. We note that as is the case in tensor-triangular geometry, $(\Spc(\catK),\supp)$ recovers the lattice of semiprime\footnote{i.e., an intersection of primes. All thick $\otimes$-ideals are semiprime if $\catK$ is rigid.} thick $\otimes$-ideals of $\catK$, see \cite[Theorem A]{DDM25} for details.
    \end{definition}

    \subsection{The Heisenberg category}

    Now, we review the construction of the Heisenberg category as defined in \cite{Kho14} and the diagrammatics which are necessary for this paper, following the streamlined presentation of \cite[Remark 1.5(2)]{Bru18}. We describe $\Heis := \Heis_{-1}(0)$ in the notation of \textit{loc. cit.} By ``Heisenberg category'' we refer to the category $\Heis$, as opposed to its Karoubi envelope, as in \cite{Kho14}.

    \begin{definition}
        Fix a commutative ring $\bbF$. We define the additive, $\bbF$-linear, strict monoidal category $\Heis$ generated by two objects $\up$ and $\down$. That is, an object of $\Heis$ is a finite direct sum of tensor products $X = X_1\otimes \cdots \otimes X_n$, where $X$ is a finite word with each character $X_i \in \{\up,\down\}$. The value of $X_i$ is its \textit{orientation}, each $X_i$ is either upwards or downwards oriented.

        The morphisms of $\Heis$ are generated by the identity morphisms on $\up$ and $\down$,
        \[
            \id_\up =
            \mathord{
                \begin{tikzpicture}[baseline = 0]
                    \draw[->,thick,darkblue] (0.08,-.3) to (0.08,.4);
                \end{tikzpicture}
            }
            \;,
            \qquad\quad \id_\down =
            \mathord{
                \begin{tikzpicture}[baseline = 0]
                	\draw[<-,thick,darkblue] (0.08,-.3) to (0.08,.4);
                \end{tikzpicture}
            }\:,
        \]
        and morphisms $s: \up\otimes\up \to \up\otimes\up$, $c: \bbone \to \down\otimes\up$, $d: \up\otimes \down\to \bbone$, $c': \bbone\to \up\otimes\down$, and $d': \down\otimes\up \to \bbone$,

        \[
            s =
            \mathord{
                \begin{tikzpicture}[baseline = 0]
                    \draw[->,thick,darkblue] (0.28,-.3) to (-0.28,.4);
                    \draw[thick,darkblue,->] (-0.28,-.3) to (0.28,.4);
                \end{tikzpicture}
            }
            \:,
            \qquad c =
            \mathord{
                \begin{tikzpicture}[baseline = 1mm]
                    \draw[<-,thick,darkblue] (0.4,0.4) to[out=-90, in=0] (0.1,0);
                    \draw[-,thick,darkblue] (0.1,0) to[out = 180, in = -90] (-0.2,0.4);
                \end{tikzpicture}
            }
            \:,
            \qquad
            d =
            \mathord{
                \begin{tikzpicture}[baseline = 1mm]
                    \draw[<-,thick,darkblue] (0.4,0) to[out=90, in=0] (0.1,0.4);
                    \draw[-,thick,darkblue] (0.1,0.4) to[out = 180, in = 90] (-0.2,0);
                \end{tikzpicture}
            }
            \:,
            \qquad c' =
            \mathord{
                \begin{tikzpicture}[baseline = 1mm]
                    \draw[-,thick,darkblue] (0.4,0.4) to[out=-90, in=0] (0.1,0);
                    \draw[->,thick,darkblue] (0.1,0) to[out = 180, in = -90] (-0.2,0.4);
                \end{tikzpicture}
            }
            \:,
            \qquad d' =
            \mathord{
                \begin{tikzpicture}[baseline = 1mm]
                    \draw[-,thick,darkblue] (0.4,0) to[out=90, in=0] (0.1,0.4);
                    \draw[->,thick,darkblue] (0.1,0.4) to[out = 180, in = 90] (-0.2,0);
                \end{tikzpicture}
            }
            \:.
        \]

        These diagrams record the relations that the morphisms satisfy. Tensor products, i.e., \textit{horizontal composition} $a \otimes b$, of two morphisms is $a$ drawn to the left of $b$, and \textit{vertical composition} $a \circ b$ is $a$ drawn above $b$ (when this makes sense). We define the morphisms $t: \up\otimes\down \to \down\otimes\up$ and $t': \down\otimes\up\to \up\otimes\down$ by
        \[
            t=\mathord{
                \begin{tikzpicture}[baseline = 0]
                	\draw[<-,thick,darkblue] (0.28,-.3) to (-0.28,.4);
                	\draw[->,thick,darkblue] (-0.28,-.3) to (0.28,.4);
                \end{tikzpicture}
            }
            :=
            \mathord{
                \begin{tikzpicture}[baseline = 0]
                	\draw[->,thick,darkblue] (0.3,-.5) to (-0.3,.5);
                	\draw[-,thick,darkblue] (-0.2,-.2) to (0.2,.3);
                        \draw[-,thick,darkblue] (0.2,.3) to[out=50,in=180] (0.5,.5);
                        \draw[->,thick,darkblue] (0.5,.5) to[out=0,in=90] (0.8,-.5);
                        \draw[-,thick,darkblue] (-0.2,-.2) to[out=230,in=0] (-0.5,-.5);
                        \draw[-,thick,darkblue] (-0.5,-.5) to[out=180,in=-90] (-0.8,.5);
                \end{tikzpicture}
            }\,,
            \qquad\quad
            t' = \mathord{
            \begin{tikzpicture}[baseline = 0]
            	\draw[->,thick,darkblue] (0.28,-.3) to (-0.28,.4);
            	\draw[<-,thick,darkblue] (-0.28,-.3) to (0.28,.4);
            \end{tikzpicture}
            }
            :=
            \mathord{
                \begin{tikzpicture}[baseline = 0, xscale = -1]
                	\draw[->,thick,darkblue] (0.3,-.5) to (-0.3,.5);
                	\draw[-,thick,darkblue] (-0.2,-.2) to (0.2,.3);
                        \draw[-,thick,darkblue] (0.2,.3) to[out=50,in=180] (0.5,.5);
                        \draw[->,thick,darkblue] (0.5,.5) to[out=0,in=90] (0.8,-.5);
                        \draw[-,thick,darkblue] (-0.2,-.2) to[out=230,in=0] (-0.5,-.5);
                        \draw[-,thick,darkblue] (-0.5,-.5) to[out=180,in=-90] (-0.8,.5);
                \end{tikzpicture}
            }\,.
        \]
        Then we impose the following \textit{degenerate Hecke relations},

        \[
            \mathord{
                \begin{tikzpicture}[baseline = -1mm]
                	\draw[->,thick,darkblue] (0.28,0) to[out=90,in=-90] (-0.28,.6);
                	\draw[->,thick,darkblue] (-0.28,0) to[out=90,in=-90] (0.28,.6);
                	\draw[-,thick,darkblue] (0.28,-.6) to[out=90,in=-90] (-0.28,0);
                	\draw[-,thick,darkblue] (-0.28,-.6) to[out=90,in=-90] (0.28,0);
                \end{tikzpicture}
            }=
            \mathord{
                \begin{tikzpicture}[baseline = -1mm]
                	\draw[->,thick,darkblue] (0.18,-.6) to (0.18,.6);
                	\draw[->,thick,darkblue] (-0.18,-.6) to (-0.18,.6);
                \end{tikzpicture}
            }\,,\qquad
            \mathord{
                \begin{tikzpicture}[baseline = -1mm]
                	\draw[<-,thick,darkblue] (0.45,.6) to (-0.45,-.6);
                	\draw[->,thick,darkblue] (0.45,-.6) to (-0.45,.6);
                        \draw[-,thick,darkblue] (0,-.6) to[out=90,in=-90] (-.45,0);
                        \draw[->,thick,darkblue] (-0.45,0) to[out=90,in=-90] (0,0.6);
                \end{tikzpicture}
            }
            =
            \mathord{
                \begin{tikzpicture}[baseline = -1mm]
                	\draw[<-,thick,darkblue] (0.45,.6) to (-0.45,-.6);
                	\draw[->,thick,darkblue] (0.45,-.6) to (-0.45,.6);
                        \draw[-,thick,darkblue] (0,-.6) to[out=90,in=-90] (.45,0);
                        \draw[->,thick,darkblue] (0.45,0) to[out=90,in=-90] (0,0.6);
                \end{tikzpicture}
            },\,
        \]
        and the \textit{right adjunction relations},\footnote{One also has left adjunction relations, which are obtained by flipping the diagrams horizontally.}
        \[
        \mathord{
            \begin{tikzpicture}[baseline = 0]
              \draw[->,thick,darkblue] (0.3,0) to (0.3,.4);
            	\draw[-,thick,darkblue] (0.3,0) to[out=-90, in=0] (0.1,-0.4);
            	\draw[-,thick,darkblue] (0.1,-0.4) to[out = 180, in = -90] (-0.1,0);
            	\draw[-,thick,darkblue] (-0.1,0) to[out=90, in=0] (-0.3,0.4);
            	\draw[-,thick,darkblue] (-0.3,0.4) to[out = 180, in =90] (-0.5,0);
              \draw[-,thick,darkblue] (-0.5,0) to (-0.5,-.4);
            \end{tikzpicture}
            }
            =
            \mathord{
            \begin{tikzpicture}[baseline=0]
              \draw[->,thick,darkblue] (0,-0.4) to (0,.4);
            \end{tikzpicture}
            }\,,\qquad
            \mathord{
                \begin{tikzpicture}[baseline = 0]
                  \draw[->,thick,darkblue] (0.3,0) to (0.3,-.4);
                	\draw[-,thick,darkblue] (0.3,0) to[out=90, in=0] (0.1,0.4);
                	\draw[-,thick,darkblue] (0.1,0.4) to[out = 180, in = 90] (-0.1,0);
                	\draw[-,thick,darkblue] (-0.1,0) to[out=-90, in=0] (-0.3,-0.4);
                	\draw[-,thick,darkblue] (-0.3,-0.4) to[out = 180, in =-90] (-0.5,0);
                  \draw[-,thick,darkblue] (-0.5,0) to (-0.5,.4);
                \end{tikzpicture}
            }
            =
            \mathord{
                \begin{tikzpicture}[baseline=0]
                  \draw[<-,thick,darkblue] (0,-0.4) to (0,.4);
                \end{tikzpicture}
            }\,.
        \]
        Finally we assert four additional relations:
        \[
            \mathord{
            \begin{tikzpicture}[baseline = 0mm]
            	\draw[->,thick,darkblue] (0.28,0) to[out=90,in=-90] (-0.28,.6);
            	\draw[-,thick,darkblue] (-0.28,0) to[out=90,in=-90] (0.28,.6);
            	\draw[<-,thick,darkblue] (0.28,-.6) to[out=90,in=-90] (-0.28,0);
            	\draw[-,thick,darkblue] (-0.28,-.6) to[out=90,in=-90] (0.28,0);
            \end{tikzpicture}
            }
            =\mathord{
            \begin{tikzpicture}[baseline = 0]
            	\draw[<-,thick,darkblue] (0.08,-.6) to (0.08,.6);
            	\draw[->,thick,darkblue] (-0.28,-.6) to (-0.28,.6);
            \end{tikzpicture}
            }\,,\qquad
            \mathord{
            \begin{tikzpicture}[baseline = 0mm]
            	\draw[-,thick,darkblue] (0.28,0) to[out=90,in=-90] (-0.28,.6);
            	\draw[->,thick,darkblue] (-0.28,0) to[out=90,in=-90] (0.28,.6);
            	\draw[-,thick,darkblue] (0.28,-.6) to[out=90,in=-90] (-0.28,0);
            	\draw[<-,thick,darkblue] (-0.28,-.6) to[out=90,in=-90] (0.28,0);
            \end{tikzpicture}
            }
            =\mathord{
            \begin{tikzpicture}[baseline = 0]
            	\draw[->,thick,darkblue] (0.08,-.6) to (0.08,.6);
            	\draw[<-,thick,darkblue] (-0.28,-.6) to (-0.28,.6);
            \end{tikzpicture}
            }
            -\mathord{
            \begin{tikzpicture}[baseline=-.5mm]
            	\draw[<-,thick,darkblue] (0.3,0.6) to[out=-90, in=0] (0,.1);
            	\draw[-,thick,darkblue] (0,.1) to[out = 180, in = -90] (-0.3,0.6);
            	\draw[-,thick,darkblue] (0.3,-.6) to[out=90, in=0] (0,-0.1);
            	\draw[->,thick,darkblue] (0,-0.1) to[out = 180, in = 90] (-0.3,-.6);
            \end{tikzpicture}}\,,\qquad
            \mathord{
            \begin{tikzpicture}[baseline = -0.5mm]
            	\draw[<-,thick,darkblue] (0,0.6) to (0,0.3);
            	\draw[-,thick,darkblue] (0,0.3) to [out=-90,in=0] (-.3,-0.2);
            	\draw[-,thick,darkblue] (-0.3,-0.2) to [out=180,in=-90](-.5,0);
            	\draw[-,thick,darkblue] (-0.5,0) to [out=90,in=180](-.3,0.2);
            	\draw[-,thick,darkblue] (-0.3,.2) to [out=0,in=90](0,-0.3);
            	\draw[-,thick,darkblue] (0,-0.3) to (0,-0.6);
            \end{tikzpicture}
            }=
            0,\qquad
            \mathord{
            \begin{tikzpicture}[baseline = 1.25mm]
              \draw[->,thick,darkblue] (0.2,0.2) to[out=90,in=0] (0,.4);
              \draw[-,thick,darkblue] (0,0.4) to[out=180,in=90] (-.2,0.2);
            \draw[-,thick,darkblue] (-.2,0.2) to[out=-90,in=180] (0,0);
              \draw[-,thick,darkblue] (0,0) to[out=0,in=-90] (0.2,0.2);
            \end{tikzpicture}
            }= \id_\bbone.
        \]
    \end{definition}

    \begin{proposition}
        The \textit{Heisenberg relation} $ \up\otimes\down \oplus \bbone \cong \down\otimes \up$ holds.
    \end{proposition}
    \begin{proof}
        Indeed, the isomorphism is given by \[
        \left[\:
        \mathord{
        \begin{tikzpicture}[baseline = 0]
        	\draw[<-,thick,darkblue] (0.28,-.3) to (-0.28,.4);
        	\draw[->,thick,darkblue] (-0.28,-.3) to (0.28,.4);
        \end{tikzpicture}
        }\:\:\:
        \mathord{
        \begin{tikzpicture}[baseline = -0.9mm]
        	\draw[<-,thick,darkblue] (0.4,0.2) to[out=-90, in=0] (0.1,-.2);
        	\draw[-,thick,darkblue] (0.1,-.2) to[out = 180, in = -90] (-0.2,0.2);
        \end{tikzpicture}
        }
        \:
        \right]
        :\up \otimes \down \oplus
        \bbone
        \xrightarrow{}
         \down \otimes  \up
        \]
        with inverse
        \[
        \left[\:
        \mathord{
        \begin{tikzpicture}[baseline = 0]
        	\draw[<-,thick,darkblue] (0.28,-.3) to (-0.28,.4);
        	\draw[->,thick,darkblue] (-0.28,-.3) to (0.28,.4);
        \end{tikzpicture}
        }\:\:\:
        \mathord{
        \begin{tikzpicture}[baseline = -0.9mm]
        	\draw[<-,thick,darkblue] (0.4,0.2) to[out=-90, in=0] (0.1,-.2);
        	\draw[-,thick,darkblue] (0.1,-.2) to[out = 180, in = -90] (-0.2,0.2);
        \end{tikzpicture}
        }
        \:
        \right]
        ^{-1} :=
        \left[\begin{array}{l}
        \displaystyle\mathord{
        \begin{tikzpicture}[baseline = 0]
        	\draw[->,thick,darkblue] (0.28,-.3) to (-0.28,.4);
        	\draw[<-,thick,darkblue] (-0.28,-.3) to (0.28,.4);
        \end{tikzpicture}
        }\\\\
        \mathord{
        \begin{tikzpicture}[baseline = 1mm]
        	\draw[-,thick,darkblue] (0.4,0) to[out=90, in=0] (0.1,0.4);
        	\draw[->,thick,darkblue] (0.1,0.4) to[out = 180, in = 90] (-0.2,0);
        \end{tikzpicture}
        }
        \end{array}\right]: \down\otimes\up\xrightarrow{}\up\otimes\down\oplus\bbone.
        \]

    \end{proof}

    We set \[x = \mathord{
        \begin{tikzpicture}[baseline = -.5mm]
        	\draw[->,thick,darkblue] (0.08,-.4) to (0.08,.5);
            \node at (0.08,0.05) {$\dot$};
        \end{tikzpicture}
        }
        :=
        \mathord{
            \begin{tikzpicture}[baseline = -0.5mm, xscale = -1]
            	\draw[<-,thick,darkblue] (0,0.6) to (0,0.3);
            	\draw[-,thick,darkblue] (0,0.3) to [out=-90,in=0] (-.3,-0.2);
            	\draw[-,thick,darkblue] (-0.3,-0.2) to [out=180,in=-90](-.5,0);
            	\draw[-,thick,darkblue] (-0.5,0) to [out=90,in=180](-.3,0.2);
            	\draw[-,thick,darkblue] (-0.3,.2) to [out=0,in=90](0,-0.3);
            	\draw[-,thick,darkblue] (0,-0.3) to (0,-0.6);
            \end{tikzpicture}
        }\,,
    \]
    and denote the $n$th power $x^{\circ n}$ of $x$ under vertical composition diagrammatically by labeling the dot with the multiplicity $n$. We have the following counterparts to $x$ and $s$, denoted $x': \down\to \down$ and $s': \down\otimes\down\to \down\otimes\down$ respectively via
    \[
        x' = \mathord{
            \begin{tikzpicture}[baseline = 0]
        	\draw[<-,thick,darkblue] (0.08,-.3) to (0.08,.4);
              \node at (0.08,0.1) {$\dot$};
        \end{tikzpicture}
        }
        :=
        \mathord{
        \begin{tikzpicture}[baseline = 0]
          \draw[->,thick,darkblue] (0.3,0) to (0.3,-.4);
        	\draw[-,thick,darkblue] (0.3,0) to[out=90, in=0] (0.1,0.4);
        	\draw[-,thick,darkblue] (0.1,0.4) to[out = 180, in = 90] (-0.1,0);
        	\draw[-,thick,darkblue] (-0.1,0) to[out=-90, in=0] (-0.3,-0.4);
        	\draw[-,thick,darkblue] (-0.3,-0.4) to[out = 180, in =-90] (-0.5,0);
          \draw[-,thick,darkblue] (-0.5,0) to (-0.5,.4);
           \node at (-0.1,0) {$\dot$};
        \end{tikzpicture}
        }
        =
        \mathord{
        \begin{tikzpicture}[baseline = 0]
          \draw[-,thick,darkblue] (0.3,0) to (0.3,.4);
        	\draw[-,thick,darkblue] (0.3,0) to[out=-90, in=0] (0.1,-0.4);
        	\draw[-,thick,darkblue] (0.1,-0.4) to[out = 180, in = -90] (-0.1,0);
        	\draw[-,thick,darkblue] (-0.1,0) to[out=90, in=0] (-0.3,0.4);
        	\draw[-,thick,darkblue] (-0.3,0.4) to[out = 180, in =90] (-0.5,0);
          \draw[->,thick,darkblue] (-0.5,0) to (-0.5,-.4);
           \node at (-0.1,0) {$\dot$};
        \end{tikzpicture}
        }\:,\qquad
        s' = \mathord{
        \begin{tikzpicture}[baseline = 0]
        	\draw[<-,thick,darkblue] (0.28,-.3) to (-0.28,.4);
        	\draw[<-,thick,darkblue] (-0.28,-.3) to (0.28,.4);
        \end{tikzpicture}
        }:=
        \mathord{
        \begin{tikzpicture}[baseline = 0]
        	\draw[<-,thick,darkblue] (0.3,-.5) to (-0.3,.5);
        	\draw[-,thick,darkblue] (-0.2,-.2) to (0.2,.3);
                \draw[-,thick,darkblue] (0.2,.3) to[out=50,in=180] (0.5,.5);
                \draw[->,thick,darkblue] (0.5,.5) to[out=0,in=90] (0.8,-.5);
                \draw[-,thick,darkblue] (-0.2,-.2) to[out=230,in=0] (-0.5,-.5);
                \draw[-,thick,darkblue] (-0.5,-.5) to[out=180,in=-90] (-0.8,.5);
        \end{tikzpicture}
        }=
        \mathord{
        \begin{tikzpicture}[baseline = 0]
        	\draw[->,thick,darkblue] (0.3,.5) to (-0.3,-.5);
        	\draw[-,thick,darkblue] (-0.2,.2) to (0.2,-.3);
                \draw[-,thick,darkblue] (0.2,-.3) to[out=130,in=180] (0.5,-.5);
                \draw[-,thick,darkblue] (0.5,-.5) to[out=0,in=270] (0.8,.5);
                \draw[-,thick,darkblue] (-0.2,.2) to[out=130,in=0] (-0.5,.5);
                \draw[->,thick,darkblue] (-0.5,.5) to[out=180,in=-270] (-0.8,-.5);
        \end{tikzpicture}
        }\:.
    \]

    We introduce some terminology which will be of import in the sequel. For the sake of simplicity, in the following definition when we say ``connected,'' we discount string crossings.

    \begin{definition}
        If $S$ is a connected string (either with endpoints at the boundaries or a closed loop) in a diagram, we say $S$ is:
        \begin{itemize}
            \item a \textit{cup} if $S$ has both endpoints on the upper boundary (e.g., $c,c'$);
            \item a \textit{cap} if $S$ has both endpoints on the lower boundary (e.g., $d,d'$);
            \item a \textit{bridge} if $S$ has endpoints on both lower and upper boundaries (e.g. $\id_\up, x$);
            \item a \textit{bubble} if it has no endpoints (e.g., $d\circ c'$, $c\circ d' = \id_\bbone$).
        \end{itemize}

        Any connected strand must be one of these four types. A diagram without endpoints is an endomorphism of the tensor unit $\bbone$, in which case, every string is a bubble. The calculus of \textit{bubble moves}, i.e., crossing bubbles over strings, is nontrivial and follows from other relations, but we omit its review, as it will not be necessary for this paper.

        For cups, caps, and bubbles, one can speak of their \textit{orientation} as being \textit{clockwise} or \textit{counter-clockwise}, and for bridges, one can speak of their orientation as \textit{upwards} or \textit{downwards}. For strings with endpoints, the orientations can be read purely off of the endpoint orientation: the orientation of a bridge matches the orientation of both of its endpoints, and a clockwise-oriented cup or cap has its up-oriented endpoint to the left of its down-oriented endpoint.

        Similarly, there are two possible types of \textit{curls} on strings: a clockwise curl (e.g., $x$, $x'$) and a counterclockwise curl. Any diagram containing any left curl subdiagram is necessarily zero by the relations.
        A string with $n$ dots can be closed into either a clockwise- or counterclockwise-oriented circle with $n$ dots. However a counterclockwise loop with a dot, i.e. a \textit{figure-eight diagram}, is the 0 morphism (thus any diagram containing such a loop is also 0), and a counterclockwise loop with 2 or more dots can be expressed as a linear combination of products of clockwise loops, see \cite{Kho14} for details.

        Given an unreduced composition of diagrams $g \circ f$ and a continuous strand $S$ in $g\circ f$, it will be important to analyze the sub-strands of $S$ in $g$ and $f$ that comprise $S$. We denote the \textit{set of connected substrands} of $S$ in $g$ and $f$ respectively as $S \cap g$ and $S \cap f$ (note that $S \cap g$ or $S\cap f$ need not be connected). This notation depends on how $g\circ f$ is drawn (for instance, applying the right adjunction relation may remove cups or caps), we assume $g \circ f$ is unreduced.
    \end{definition}

    \subsection{The standard basis}

    We now review an $\bbF$-basis of $\Hom_{\Heis}(X, Y)$ for any finite words $X,Y$ with characters in $\up$ and $\down$, following \cite{Kho14}. An $(X,Y)$-matching is a bijection \[\{i \mid X_i = \up\}\sqcup \{j \mid Y_j = \down\} \xrightarrow{\cong} \{i \mid X_i = \down\} \sqcup \{j \mid Y_j = \up\}.\] If no such bijection exists, $\Hom_{\Heis}(X, Y) = 0$.

    A \textit{reduced lift} of an $(X,Y)$-matching is a diagram representing a morphism $X \to {Y}$ such that
    \begin{itemize}
        \item the endpoints of each string in the diagram are paired under the matching;
        \item any two strings intersect at most once;
        \item there are no self-intersections (e.g., no dots/loops);
        \item there are no triple intersections;
        \item there are no bubbles.
    \end{itemize}

    Let $B_0(X,Y)$ be a set consisting of a reduced lift for each of the $(X,Y)$-matchings. Note the vertical composition of two reduced lifts need not be a reduced lift. For each element of $B_0(X,Y)$, pick a distinguished point on each string that is away from crossings and critical points, and near the out endpoint of the string. For each distinguished point, assign a dot with a nonnegative integer to it. Then, choose a finite sub-multiset of $\bbN_{\geq 0}$ and for each element, draw a clockwise bubble with a dot labeled with the corresponding nonnegative integer in the rightmost portion of the diagram, with no intersections. The resulting set of diagrams is denoted $B(X,Y)$.
    If $n$ denotes the total number of $\up$'s in $X$ and $\down$'s in $Y$, then $B(X, Y)$ is parametrized by the $(X,Y)$-matchings (of which there are $n!$), a sequence of $n$ nonnegative integers describing the number of dots on each string, and another finite sequence of nonnegative integers listing the number of clockwise oriented bubbles with no dots, one dot, and so on. In particular, $B_0(X,Y)$ is the subset of $B(X,Y)$ which are only parametrized by the $(X,Y)$-matchings.

    \begin{proposition}{\cite[Proposition 5]{Kho14}}
        For any sequences $X,Y$, the set $B(X,Y)$ constitutes a $\bbF$-basis of the $\bbF$-module $\Hom_{\Heis}(X,Y).$
    \end{proposition}

    \section{The main results}

    We first describe the semisimple subcategory of $\Heis$. To do so, we describe for all sequences $X, Y$, a subset of $B_0(X,Y)$ compatible under vertical and horizontal composition. We first record a few observations.

    \begin{lemma}\label{lem:substringobservations}
        Let $f\in B(X,Y)$ and $g \in B(Y,Z)$ be composable diagrams, and let $S$ be a string in $g\circ f$.
        \begin{enumerate}
            \item If $S$ is a non-self-intersecting clockwise (resp. counterclockwise)-oriented cup, cap, or bubble, then there exists a substring $S'\subset S$ in $f$ or $g$ which is a cup, cap, or bubble with the same orientation.
            \item If $S$ is a bridge, then there exist unique substrings $S' \in S \cap f$ and $S'' \in S \cap g$ which are also bridges. Moreover, $S'$ and $S''$ have the same orientation as $S$.
        \end{enumerate}
    \end{lemma}
    \begin{proof}
        For (a), let $S$ be a cup in $g\circ f$. If $S \subset g$, there is nothing to prove, so assume $S$ is partially contained in $f$, in which case $S \cap g$ consists of two bridges and possibly some caps, and $S \cap f$ consists of at least one cup. Suppose $g\colon Y \to Z$ with $Y = Y_1 \otimes \cdots \otimes Y_m$, $Z = Z_1 \otimes \cdots \otimes Z_n$ with the two bridges in $S \cap g$ connecting $Y_i$, $Z_{i'}$ and $Y_j$, $Z_{j'}$ respectively. Assume without loss of generality that $i < j$, then $i' < j'$, as $S$ is non-intersecting. Moreover, we have orientation matchings $Y_i=Z_{i'}$ and $Y_j=Z_{j'}$, $Y_i \neq Y_j$ are opposites, and the choice of orientation of $Y_i$ (or any of the other letters) determines the orientation of $S$.

        Now, because $S$ is non-intersecting, there must exist a cup substring $S'$ in $S \cap f$ connecting some nodes $Y_k$ and $Y_l$, with $k \leq i$ and $l > i$. Indeed, if the cup connecting $Y_i$ in $f$ has other exit node further to the right, then setting $k = i$ suffices. Otherwise, eventually $S$ must cross to the right of $Y_i$ (as $S$ must reach $Y_j$), and must do so in $f$, since otherwise $S$ would intersect itself at the bridge with exit nodes $Y_i$ and $Z_{i'}$. Take the first time, moving in the direction on $S$ away from $Y_i$, that it does so: suppose it happens at a cup  $S'$ in $S \cap f$ passing through $Y_k$ with $k < i$ and $Y_l$ with $l > i$. By a parity argument, $Y_k$ has the same orientation as $Y_i$, and because $k < l$, we conclude $S'$ as a cup in $f$ has the same orientation as $S$. The analogous argument for $S$ a cap holds.

        Finally, for $S$ a bubble, if $S \subseteq f$ or $S \subseteq g$ there is nothing to show, so suppose $S \cap f\neq \emptyset$ and $S \cap g \neq \emptyset$. Assuming $g \colon Y \to Z$ as before, take the least $i \in \bbN$ such that $S$ passes through $Y_i$. Note the orientation of $Y_i$ determines the orientation of $S$; then either the corresponding cup in $S\cap f$ passing through $Y_i$ or cup in $S\cap g$ have the same orientation as $S$. Thus (a) is shown.

        For (b), the corresponding bridges in $f$ and $g$ are the first and last occurring maximal substrings in $S$, it is evident that they must have the same orientation as $S$, and they are unique as all other substrings cannot reach the outer boundaries of $f$ and $g$, hence must be cups and caps.
    \end{proof}

    Set $A_1 :=\bbZ[x,\partial]/\langle \partial x - x\partial  - 1\rangle$, the \textit{integral Weyl algebra on one variable}.

    \begin{theorem}\label{thm:semisimplesubcat}
        There exists a replete semisimple monoidal subcategory $\Heis_{s} \subseteq \Heis$ given as follows. The objects of $\Heis_{s}$ are the objects of $\calH$ and the morphism spaces $\Hom(X, Y)$ has a $\bbF$-basis described by the subset of $B_0(X, Y)\subset B(X,Y)$ consisting of diagrams for which the following hold:
        \begin{enumerate}
            \item There are no dots or bubbles. That is, the diagram belongs to $B_0(X,Y)$;
            \item There are no clockwise-oriented cups or caps (e.g., $c',d$ are excluded). In particular, all cups and caps correspond to a pair $\up,\down$ for which the $\down$ occurs on the left of $\up$;
            \item Bridges may cross only if they have opposite orientation (e.g., $s, s'$ are excluded, $t, t'$ are allowed).
        \end{enumerate}
        The Heisenberg isomorphism $\down \otimes \up\cong \up \otimes \down\oplus \bbone$ holds in $\Heis_{s}$. Moreover, the simple objects are given precisely by $\up^{\otimes i}\otimes \down^{\otimes j}$ for nonnegative integers $i,j \geq 0$. Finally, we have an isomorphism $K_0(\Heis_s) \cong A_1$.
    \end{theorem}

    \begin{proof}
        We must verify that $\Heis_s$ is closed under vertical and horizontal composition. Horizontal composition is straightforward. We verify (a)-(c) for vertical composition, in particular verifying that if two diagrams $g, f$ compose to something breaking one of the conditions, the composition is equivalent to the 0 morphism.

        For condition (a), let $f\in B_0(X, Y)$ and $g \in B_0(Y, Z)$ be composable morphisms in $\Heis_s$ satisfying (a)-(c). If $g\circ f$ contains a non-self-intersecting bubble $S$, then since neither $f$ nor $g$ have any bubbles, $S \cap f$ and $S \cap G$ must consist of cups and caps, hence by \Cref{lem:substringobservations}(a), $S$ must be counterclockwise-oriented, and therefore can be removed. Now, suppose $g\circ f$ contains a self-intersecting strand $S$, intersecting itself any number of times. Then there exists a connected substring $S'\subseteq S$ which intersects itself exactly once and contains a closed loop; it is locally a curl (note the choice of intersection need not be unique if $S$ is a bubble, but is unique otherwise). For instance, in the example below, one chooses the blue substrand. Because neither $f$ nor $g$ contain self-intersections, $S'$ cannot be contained entirely in $f$ or $g$. Suppose without loss of generality that the intersection in $S'$ happens in $f$, then $S' \cap g$ contains at least one cup, and since $S' \cap g \subseteq g$, the cup is necessarily counterclockwise-oriented. Therefore, $S'$ is a counterclockwise-oriented loop, hence $g\circ f = 0$. Therefore there are no dots or self-intersecting bubbles in $g\circ f$ after reduction.

        For condition (b) let $f\in B_0(X, Y)$ and $g \in B_0(Y, Z)$ be as before, and suppose $S$ in $g \circ f$ is a cup. If $S\in g$, there is nothing to show, so suppose $S \cap f \neq 0$. If $S$ is non-self-intersecting, then it follows from \Cref{lem:substringobservations}(a) that $S$ is counterclockwise oriented, since all caps and cups in $f$ and $g$ are. Otherwise, suppose $S$ is self-intersecting. In particular, $S$ contains a curl, and by condition (a), which we have shown is closed under composition, the self-intersection is a counterclockwise curl, and $g\circ f = 0$. The cap case follows analogously.

        For condition (c), let $f\in B_0(X, Y)$ and $g \in B_0(Y, Z)$ be as before, and let $S$ be a cup in $f$. If $S$ belongs to a bubble in $g\circ f$, there is nothing to show, since if the bubble is non-self-intersecting it may be removed, and if it is self-intersecting $g\circ f = 0$ by (a) and (b). If $S$ meets two noncrossing bridges, then in $g\circ f$, the string containing $S$ remains a cup, and if $S$ meets two bridges which cross once (the maximum number of times they can cross, since $g \in B_0(Y, Z)$, then $S$ belongs to a left loop in $g\circ f$, hence $g\circ f = 0$.
        Finally, if $S$ meets a cap at one endpoint but is not a , then we claim $g \circ f = 0$ regardless of what $S$ meets at the other endpoint. Indeed, let $S$ meet the cap $S'$ in $g$. Assume without loss of generality that $S$ and $S'$ meet at a $\down$ endpoint, the $\up$ case follows analogously. In this case, the other $\up$ endpoints of $S$ and $S'$ both lay to the right of the $\down$ at which $S, S'$ meet, since both $S$ and $S'$ are counterclockwise-oriented. An example is provided below, with $S$ the brown substring and $S'$ the purple substring.

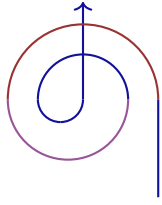
        \[
            \mathord{
            \begin{tikzpicture}[baseline = 0]
                \draw[->,thick,darkblue] (1,0) to (1,1.3);
                \draw[-,thick,darkblue] (1,0) to[out = 270, in = 0] (.7,-.3);
                \draw[-,thick,darkblue] (.7,-.3) to[out=180, in = 270] (.4, 0);
                \draw[-,thick,darkblue] (.4,0) to[out=90, in = 180] (1,.6);
                \draw[-,thick,darkblue] (1,.6) to[out = 0, in = 90] (1.6,0);
                \draw[-,thick,darkpurple] (1.6,0) to[out = 270, in = 0] (.8,-.8);
                \draw[-,thick,darkpurple] (.8,-.8) to[out = 180, in = 270] (0,0);
                \draw[-,thick,darkred] (0,0) to[out = 90, in = 180] (1,1);
                \draw[-,thick,darkred] (1,1) to[out = 0, in =90] (2,0);
                \draw[-,thick,darkblue] (2,0) to (2,-1.3);
            \end{tikzpicture}
            }
        \]

        It follows by parity of the remaining available nodes enclosed inside of $S$ and $S'$ that there must be another string connected to $S$ and $S'$ in $g\circ f$ crossing one of $S$ and $S'$, as in the above example. Indeed, the string containing $S$ and $S'$ must either terminate above $S$, below $S'$, or contain no endpoints, in which case the string meets $S$ again. In any of these cases, the crossing forms a counterclockwise curl, and thus $g \circ f = 0$.

        The analogous argument as above may be performed for $S$ a cap in $g$. By ruling out the above cases, we conclude that if there exist crossing bridges $B$, $B'$ in $g \circ f$ with $g\circ f$ nonzero, then the corresponding strings in $B\cap g$, $B\cap f$ and $B'\cap g$, $B'\cap f$ must all be bridges as well. However, by \Cref{lem:substringobservations}(b), a bridge $S \in f$ cannot belong to a bridge with opposite orientation in the composition $g\circ f$. Since the crossing bridges of $g \circ f$ must cross in either $g$ or $f$ (or both, in which case they can be unwound using the local relations), and those must have opposite orientation, the crossing in $g \circ f$ has opposite orientation as well. Thus we have shown (c) holds under vertical composition, and conclude $\Heis_s$ is indeed a monoidal subcategory of $\Heis$.

        The Heisenberg isomorphism holds since all components of the isomorphism exist in $\Heis_s$. Let $X = \up^{\otimes i}\otimes \down^{\otimes j}$ for some non-negative $i, j \geq 0$. We show $\End(X) \cong k$, which follows from conditions (b) and (c). Since $X$ has no subsequence of the form $\down\otimes\up$, there are no caps or cups in any endomorphism of $X$ due to (b), and an inductive argument demonstrates that due to (c), the only nontrivial element of $B_0(X, X)$ belonging to $\Heis_s$ is the identity. A similar argument shows that for any other other $Y = \up^{\otimes i'}\otimes \down^{\otimes j'}$ for $(i,j)\neq (i',j')$, $\Hom(X, Y) = 0$. Conversely, the Heisenberg isomorphism demonstrates any object of $\Heis_s$ not of the form $\up^{\otimes i}\otimes \down^{\otimes j}$ decomposes, hence the simple objects are precisely those of the form $\up^{\otimes i}\otimes \down^{\otimes j}$, as desired.

        The isomorphism $A_1 \cong K_0(\Heis_s)$ now follows via the assignment $x \mapsto [\up]$ and $\partial \mapsto [\down]$. This is well-defined since the Heisenberg relation $[\down][\up] = [\up][\down] + 1$ holds, and is an isomorphism via the description of simple objects of $\Heis_s$, since $A_1$ viewed as an abelian group has a $\bbZ$-basis given by elements of the form $x^i\partial^j$ for $i,j\geq 0$.
    \end{proof}

    Note for \Cref{thm:semisimplesubcat}, we need not assume anything about our ring $\bbF$.  When $\bbF$ is a field of characteristic 0, the semisimple subcategory $\Heis_s$ contains all the isomorphism data of the Heisenberg category $\Heis$.

    \begin{proposition}\label{prop:isoclasses}
        Assume $\bbF$ is a field of characteristic 0. Let $X,Y \in \Heis$. Then $K_0(\Heis_s) \cong K_0(\Heis) \cong A_1$, and $X \cong Y \in \Heis$ if and only if $X \cong Y \in \Heis_s$.
    \end{proposition}
    \begin{proof}
        The inclusion $\Heis_s \subseteq \Heis$ induces a homomorphism on Grothendieck groups \[\phi\colon K_0(\Heis_s) \to K_0(\Heis) \subseteq K_0(\on{Kar}(\Heis)).\] If $X \cong Y \in \Heis_s$, then clearly $X \cong Y \in \Heis$ as well, so $\phi$ is surjective. On the other hand, \cite[Theorem 1.1]{BSW23} also asserts an isomorphism $A_1 \cong K_0(\Heis)$ induced by $x \mapsto [\up]$ and $\partial \mapsto [\down]$. Indeed, $K_0(\Heis)$ regarded as a subring of $K_0(\on{Kar}(\Heis))$ by the natural inclusion is isomorphic to the subring of the Heisenberg ring generated by elements $h_1^+, e_1^-$ corresponding to $\up, \down$ respectively, which has presentation $\bbZ[h_1^+,e_1^-]/\langle h_1^+e_1^- - e_1^-h_1^+ - 1\rangle$, see \cite[Remark 3.2]{BSW25}. This isomorphism is compatible with $\phi,$ thus $\phi$ is an isomorphism.
        Now, if $X \cong Y \in \Heis$, then $[X] = [Y] \in K_0(\Heis)$, therefore $[X] = [Y] \in K_0(\Heis_s)$. Because $\Heis_s$ is semisimple and Krull-Schmidt, we conclude $X \cong Y$ in $\Heis_s$, as desired.
    \end{proof}

    \begin{remark}
        We do not expect that a semisimple subcategory encoding isomorphism class data exists for $\Heis_k(\delta)$ at other central charges $k$ or nonzero parameters $\delta$. When $|k| > 1$, the Heisenberg relation requires dotted cups and caps (see \cite[1.4, 1.5]{Bru18}), and we rely on the fact that counterclockwise curls are the zero morphism, which is behavior specific to $\delta = 0$.

        The category $\Heis$ is rigid, with $X^* = \overline{X}$, where $\overline{X}$ denotes $X$ written in reverse with direction switched (e.g., $\overline{\up\otimes\down\otimes\down} = \up\otimes\up\otimes\down$). In this case, the evaluation and coevaluation morphisms are obtained by rotating the appropriate diagram 180${}^\circ$, then connecting the associated endpoints with nested cups and caps. On the other hand, $\Heis_s$ is not rigid, due to the lack of clockwise cups and caps; the only indecomposable objects that have a left dual are those of the form $\up^{\otimes i}$ and for right duals, $\down^{\otimes j}$. However, any $X$ still arises as a direct summand of ${X\otimes \overline{X}\otimes X},$ regardless of base ring $\bbF$.
    \end{remark}

    From $\Heis_s$, we obtain our example of a non-duo monoidal triangular category which still satisfies the tensor product property, thus verifying that the converse to \cite[Theorem 4.2.1]{NVY22b} does not hold.

    \begin{theorem}\label{thm:classificationofideals}
        Let $\catK := D^b(\Heis_s)$. There are 2 two-sided thick $\otimes$-ideals of $\catK$, $0$ and $\catK$, but the thick left $\otimes$-ideals are given by \[\catK \supset \langle \down\rangle_l \supset \langle \down\otimes\down\rangle_l \supset \cdots\supset 0,\] and the thick right $\otimes$-ideals are given by \[\catK \supset \langle \up \rangle_r \supset \langle \up\otimes\up \rangle_r \supset \cdots\supset 0.\]

        In particular, $\Spc(\catK) = \{\ast\}$, and for all nonzero $x \in \catK$, $\supp(x) = \Spc(\catK)$. Therefore, $\catK$ satisfies the tensor product property but is not left or right duo.
    \end{theorem}
    \begin{proof}
        We first show there are only 2 two-sided thick $\otimes$-ideals, in which case it follows that $\Spc(\catK) = \{\ast\}$. Let $\catI\subseteq \catK$ be a thick $\otimes$-ideal. From \Cref{thm:semisimplesubcat}, every object of $\catK$ is a split chain complex and therefore isomorphic as a chain complex to its homology. Therefore, it suffices to consider the indecomposable objects of $\catK$, modules regarded as chain complexes in a single degree.

        If $\catI$ is nonzero, then $X \in \catI$ with $X$ a word of length at least 1. If $X$ begins with $\up$, then $\down \otimes X \in \catI$, and therefore $X_1 \in \catI$ via the Heisenberg isomorphism, where $X_1$ denotes $X$ with the first character removed. If $X$ ends with a $\down$, then $X \otimes \up \in \catI$, hence $X^1 \in \catI$ via the Heisenberg isomorphism, where $X^1$ denotes $X$ with the last character removed. If neither case occurs, then $X$ begins with a $\down$ and ends with a $\up$ therefore, contains the subsequence $\down\otimes\up$. Therefore, one may delete the subsequence from $X$, reducing $l(X)$ by 2. In all three cases, we have deduced the existence of a $Y \in \catI$, with $l(Y)< l(X)$. Induction demonstrates that $\bbone \in \catI$, hence $\catI = \catK$ as desired. Thus $\Spc(\catK)^\vee = \ast$, so $\Spc(\catK) = \ast$ as well.

        Since $\supp(x) = \Spc(\catK)$ for every nonzero $x\in \catK$, $\catK$ trivially satisfies the tensor product property, that is, $\supp(x \otimes y) = \supp(x) \cap \supp(y)$ for all $x, y \in \catK$.

        Now we consider left and right ideals. We consider the case of right ideals, the left case follows by taking the ``dual'' left ideal $\overline{\catI}$ for any right ideal $\catI$. First observe that for any thick right $\otimes$-ideal $\catI \subseteq \catK$ and positive integer $i > 0$, $\up^{\otimes i} \not\in \catI$ if and only if $\up^{\otimes j}\not\in \catI$ for all $0 \leq j \leq i$. This is a routine verification: from this, we obtain a strict inclusion $\langle \up^{\otimes i} \rangle_r \supset \langle \up^{\otimes (i+1)} \rangle_r$ of ideals. Second, if $\up^{\otimes i}\otimes \down^{\otimes j} \in \catI$ for some positive $i, j> 0$, then $\up^{\otimes i} \in \catI$. Indeed, we have that $\up^{\otimes i}\otimes \down^{\otimes j}\otimes \up \in \catI$, and the Heisenberg isomorphism implies $\up^{\otimes i}\down^{\otimes (j-1)} \in \catI$, and the statement follows via induction. From these, it follows that every nonzero right ideal is of the form $\langle \up^{\otimes i}\rangle_r$ for some $i \geq 0$, as $\catI = \langle \up^{\otimes i}\rangle_r$, where $i$ is the least integer for which $\up^{\otimes i} \in \catI$. The result follows.
    \end{proof}

    \begin{remark}\label{rmk:classification}
        For the stable module category $\stab(\on{Coh}(\catG))$ of the coordinate algebra of a finite group scheme $\catG$, De Deyn and the author showed thick left- and right-sided $\otimes$-ideals are classified by the cohomological support variety of $\on{Coh}(\catG)$ \cite[Corollary 7.10]{DDM25}, answering a question posed by Negron--Pevtsova \cite[Question 11.1]{NP23}. The argument is mostly formal and relies on results on crossed product categories laid out by Huang--Vashaw \cite{HV25}. The reasons for this classification are highly context-specific and should not be expected to hold in general. Indeed, both $\Spc(\catK)$ and $\Spec^h(\End_{\Heis_s}(\bbone))$ are a singleton, but $\catK$ has countably many one-sided thick $\otimes$-ideals.

        Nakano--Vashaw--Yakimov proved a criterion \cite[Theorem 7.3.1]{NVY22} for when a support datum classifies right thick $\otimes$-ideals\footnote{Analogous arguments may be carried out for left thick $\otimes$-ideals.} under additional standard hypotheses. Denote the lattice of right thick $\otimes$-ideals of $\catK:= D^b(\Heis_s)$ by $\catF$. Because $\catF$ is totally ordered, it is a spatial frame, and therefore $\pt(\catF)$ is the universal open support datum by Stone duality (see \cite[Proposition 2.2]{DDM25}, here we identify $\pt(\catF)$ with the set of meet-prime elements of $\catF$), and the open subsets of $\pt(\catF)$ are in bijection with right thick $\otimes$-ideals of $\catK$. Explicitly $\pt(\catF)$ as a set identifies with $\catF \setminus \{\catK\}$ with \textit{open} sets \[\supp(\up^{\otimes n}) = \{\langle \up^{m}\rangle_r \mid m > n\}\cup \{0\}.\] Correspondingly, the support datum classifying right thick $\otimes$-ideals (in the sense of Nakano--Vashaw--Yakimov) is the Hochster dual $\pt(\catF)^\vee$, noting that $\catF$ is in fact a coherent frame, with every element compact.

        Finally, $\pt(\catF)^\vee$ coincides with the universal \textit{quasi-support datum} $\Sp(\catK)$ described in \cite[Theorem 7.2]{Mil25}. Here, $\Sp(\catK)$ is as a set the collection of all proper right thick $\otimes$-ideals. This is again a rare occurrence: one has a canonical inclusion $\pt(\catF)^\nu\hookrightarrow \Sp(\catK)$, where $\pt(\catF)^\nu$ denotes the pseudo-Hochster dual, that is, $\pt(\catF)$ with supports of objects providing a closed basis for the topology. For a general monoidal triangulated category $\catK$, this is a homeomorphism if and only if every element of $\catF$ is meet-prime, which occurs if and only if $\catF$ is totally ordered. Indeed, if $\catF$ is not totally ordered, then for any $x,y \in \catF$ for which neither $x \leq y$ nor $y \leq x$ hold, then $x \wedge y$ is not meet-prime.

    \end{remark}

    \bibliography{bib}
    \bibliographystyle{alpha}

\end{document}